\documentclass{article}
\usepackage{amsmath,amsfonts,amssymb}
\newcommand{\R}{\ensuremath{\mathbb{R}}}

\newcommand{\g}{\ensuremath{\mathfrak{g}}}
\newcommand{\PP}{\ensuremath{\mathbb{P}}}
\newcommand{\xin}{{X^{-1}}}

\newcommand{\obs}{{\hat{X}}}
\DeclareMathOperator{\inv}{inv}
\DeclareMathOperator{\Ad}{Ad}
\DeclareMathOperator{\grad}{grad}

\DeclareMathOperator{\tr}{tr}
\DeclareMathOperator{\SE}{SE}
\DeclareMathOperator{\SO}{SO}
\newcommand{\se}{\ensuremath{\mathfrak{se}}}
\newcommand{\so}{\ensuremath{\mathfrak{so}}}
\newcommand{\ddt}{\frac{\mathrm{d}}{\mathrm{d}t}}
\newtheorem{theorem}{Theorem}
\newtheorem{lemma}[theorem]{Lemma}
\newtheorem{proposition}[theorem]{Proposition}
\newtheorem{corollary}[theorem]{Corollary}
\newtheorem{remark}[theorem]{Remark}
\newtheorem{example}[theorem]{Example}
\newtheorem{definition}[theorem]{Definition}
\newenvironment{proof}{{\bf Proof: }}{\mbox{ }\hfill$\Box$\\[1ex]}

\begin{document}

\title{Gradient-like observers for invariant dynamics on a Lie group}

\author{C.~Lageman\thanks{C. Lageman is with the Centre for Mathematics and its Applications,
        The Australian National University, Canberra ACT 0200,
        Australia,
        {\tt\footnotesize lageman@maths.anu.edu.au}. His work is supported by
        the Australian Research Council Centre of Excellence for
        Mathematics and Statistics of Complex Systems.},
      J.~Trumpf\thanks{J. Trumpf is with the Department of Information Engineering,
        The Australian National University, Canberra ACT 0200,
        Australia,
        {\tt\footnotesize Jochen.Trumpf@anu.edu.au}.} and
      R.~Mahony\thanks{R. Mahony is with the Department of Engineering,
        The Australian National University, Canberra ACT 0200,
        Australia,
        {\tt\footnotesize Robert.Mahony@anu.edu.au}.}}
    
\maketitle

\begin{abstract}
This paper proposes a design methodology for non-linear 
state observers for invariant kinematic systems posed on finite
dimensional connected Lie groups, and studies the associated
fundamental system structure. The concept of synchrony of two
dynamical systems is specialised to systems on Lie groups. For
invariant systems this leads to a general factorisation theorem of a
nonlinear observer into a synchronous (internal model) term and an
innovation term.  The synchronous term is fully specified by the
system model. We propose a design methodology for the innovation
term based on gradient-like terms derived from invariant or non-invariant cost
functions. The resulting nonlinear observers have strong (almost)
global convergence properties and examples are used to demonstrate
the relevance of the proposed approach.
\end{abstract}


\section{Introduction}

There has been a surge of interest recently in nonlinear observer
design for systems on Lie groups.  A driver for recent work is the
growing demand for highly robust state estimation algorithms
for autonomous robotic systems such as unmanned aerial, ground or
submersible vehicles.  Nonlinear observer design for such
applications offers the potential of computationally simple
state-estimation algorithms with strong robustness and global
stability guarantees; as compared to the alternative of nonlinear
filter designs (eg.~extended Kalman filters \cite{2005_Anderson} or
particle filters \cite{2001_Doucet}) that provide more information
(\emph{posterior} distributions for the state estimates) but usually
require significant computational resources and rarely have strong
global stability or robustness guarantees. Previous work in
nonlinear observer design for systems on Lie groups is often closely
linked to specific applications.  In the early nineties, Salcudean
\cite{Sal1991_TAC} proposed a nonlinear observer for the attitude
estimation of a rigid-body using the unit quaternion representation
of the special-orthogonal Lie group $SO(3)$.  This work is seminal
to a series of papers that develop nonlinear attitude observers for
rigid-body dynamics
\cite{1999_Nijmeijer,ThiSan2003_TAC,2005_MahHamPfl-C64,Bonnabel2006_acc,2006_MahHamPfl_TAC,2007_Tayebi_cdc,2007_Martin_cdc},
exploiting either the unit quaternion group structure or the
rotation matrix Lie-group structure of $SO(3)$.  The resulting
attitude observers have comparable performance to state-of-the-art
nonlinear filtering techniques \cite{Crassidis07nonlinear},
generally have much stronger global stability and robustness
properties, and are simple to implement.  The full pose estimation
problem has also attracted recent attention
\cite{VikFos_2001_CDC,RehGho2003_TAC,2007_BalHamMahTru_ECC,2007_Vasconcelos_cdc}.
In this case the underlying state space is the special Euclidean
group $SE(3)$ comprising both attitude and translation of a
rigid-body. Theoretical work in this direction is less advanced due
to the more complex algebraic and geometric structure of $SE(3)$.
Recently, several authors have made a start in developing a
theoretical foundation for observer design for general systems with
invariance properties, and in particular systems with a Lie-group
state space and the natural left invariant dynamics
\cite{2008_Bonnabel_TAC,Bonnabel2005_ch}. Early work in this
direction considered uniform invariance properties across the
system, measurements and observer design. More recent work
recognises that it may be necessary to consider different invariance
properties for the system than the measurements in order to obtain
well conditioned observers \cite{2008_Bonnabel_ifac}.  Results in
this area are very recent and there is no well established observer
design methodology for invariant systems
on a Lie group, even in the case of full state measurement.

In this paper, we study the design of non-linear observers for state
space systems where the state is evolving on a finite dimensional,
connected Lie group.  We consider the case where full measurements
are available for the system kinematics and provide a
characterisation of invariant observer structures leading to an
observer design methodology.  Traditional full state observers, for
systems evolving on vector spaces, employ a design paradigm that
goes back to the work of Kalman~\cite{kalman} and
Luenberger~\cite{luenberger}: The observer system is designed as a
combination of a copy of the system or \emph{internal model} (i.e.~a
part that can in principle replicate the observed system's
trajectory), plus an innovation term which serves to drive the
observer trajectory towards the correct system trajectory in the
presence of initialization or measurement errors. We build on the
results presented in \cite{2008_Bonnabel_TAC,Bonnabel2005_ch} to
systematically study the invariance properties of internal models
and innovation terms inherent in nonlinear observer design on Lie
groups.  We define the concept of synchrony of a plant/observer pair
of systems (cf. also \cite{1997_Nijmeijer_rep}) with Lie-group
state-space in terms of constant error evolution for a canonical
invariant error term that is defined.  This leads to a definition of
internal model and a very general concept of innovation term for
plant/observer systems on Lie groups.  The internal model structure
is shown to be a copy of the invariant plant dynamics. To design the
innovation term, we utilise gradient like driving terms derived from
algebraic cost functions on the Lie group.  Under mild assumptions,
we prove almost global exponential convergence of the resulting
observers and provide example derivations for the important
applications of attitude and pose observer design for rigid-body
kinematics. Finally, we link gradient like dynamics of the canonical
error back to a specific structure of the observer consisting of a
synchronous term (internal model) and a gradient-like innovation
term. Thus, the paper provides a coherent theory of nonlinear
observer design for systems with left (resp.~right) invariant
kinematics on a Lie group for which one has a right (resp.~left)
invariant Morse-Bott cost function and full measurements.

After the introduction, Section~\ref{sec:notation} provides an
overview of the systems considered and notation used.
Section~\ref{sec:synchrony} introduces the concept of synchrony,
defines left and right invariant errors, and proves some results
concerning the canonical nature of invariant errors.
Section~\ref{sec:innovation} defines a concept of internal model and
innovation term for systems on Lie groups.
Section~\ref{sec:gradients} goes on from the characterisation of
observers to propose a specific design methodology based on using
gradients of cost functions to define the innovation terms in the
general observer structure defined in Section~\ref{sec:innovation}.
Examples of observer designs on $SO(3)$ and $SE(3)$ are provided.
The final technical section, Section~\ref{sec:grad-like}, provides
the link from gradient error dynamics back to the structure of the
observer, confirming that the design methodology proposed in
Section~\ref{sec:gradients} is the only way to obtain the natural
gradient error dynamics associated with a known cost function.  A
short paragraph of conclusions is also provided in
Section~\ref{sec:conclusions}.

A preliminary, short version of this paper has been accepted for
publication at the MTNS 2008 \cite{2008_Lageman_mtns}.

\section{Notation and problem formulation}\label{sec:notation}

Let $G$ be a finite dimensional, connected Lie group with Lie
algebra $\g$.  Denote the identity element in $G$ by $e$, and left
and right multiplication with an element $X\in G$ by $L_X$ and
$R_X$, respectively.  The tangent space $T_X G$ of $G$ at $X$ is
represented by left or right translations of the Lie algebra, i.e.
$T_e L_X \g$ or $T_e R_X \g$.  We use the simplified notation $X v$
for vectors $T_e L_X v \in T_X G$ and $v X$ for vectors $T_e R_X v
\in T_X G$ with $v\in\g$. Furthermore, we assume that there is a
Riemannian metric $\langle\cdot,\cdot\rangle$ on $G$.  Some of the
results presented in the paper will depend on invariance properties
of the metric, however, these will be stated explicitly where
required and there are no general assumptions made regarding
invariance of the metric. We use the norm $||v||^{2}=\langle
v,v\rangle$ for $v\in\g$.

Consider a left invariant system on $G$ of the form
\begin{equation}\label{eq:sys:l}
\dot X = X u ,
\end{equation}
where $u\colon \R\rightarrow\g$ is a function termed the
input signal.  An input $u$ of system (\ref{eq:sys:l}) is {\it
admissible} if solutions of the system are unique, exist for
all time and are sufficiently smooth. The system \eqref{eq:sys:l} has an equivalent right
invariant representation
\begin{equation} \label{eq:sys:r}
\dot X = v X
\end{equation}
with $v\colon \R\rightarrow\g$ an admissible input.  The input
signals $u$ and $v$ for left and right invariant representations of
the same system are related by $v = T_X R_\xin T_e L_X u = \Ad_X u$.

This paper discusses the design of observers that fuse potentially
noisy measurements of $X$ and $u$ (resp. $v$) into an estimate of
the state $X$.

\begin{example}\label{ex:so3}
Consider the design of an observer to estimate the attitude of a
rigid-body
\cite{1999_Nijmeijer,ThiSan2003_TAC,2005_MahHamPfl-C64,Bonnabel2006_acc,2006_MahHamPfl_TAC,2007_Tayebi_cdc,2007_Martin_cdc}.
The attitude of the rigid body is an element of the special
orthogonal group $\SO(3)$, represented by real, orthogonal $3\times
3$ matrices.  The Lie algebra of $SO(3)$, denoted $\so(3)$, is the
set of real, skew-symmetric $3\times 3$ matrices.  The derivative
of a curve $R\colon \R\rightarrow\SO(3)$ coincides with its
derivative in $\R^{3\times 3}$.   The maps $T_R L_S:\
R\Omega\mapsto S R\Omega$ and $T_R R_S:\ R \Omega\mapsto R\Omega S$
are given by left and right multiplication of the matrices in $T_R
\SO(3)$ with $S$, respectively.  The tangent spaces $T_R \SO(3)$ are
identified with
$$
T_R \SO(3) \equiv \{ R \Omega \;|\; \Omega\in \so(3) \} \subset \R^{3\times 3}.
$$
The special orthogonal group has a bi-invariant Riemannian metric
induced by the Euclidean metric on a the skew symmetric matrices,
i.e. $\langle R \Omega, R \Pi\rangle=\tr (\Omega^\top \Pi)$.

Consider the left invariant dynamics
$$
\dot R = R \Omega
$$
on $\SO(3)$ with $\Omega\colon \R\rightarrow \so(3)$ admissible.
This system models the kinematics of the attitude $R$ of a
coordinate frame fixed to a rigid body in 3D-space relative to an
inertial frame. Here, $\Omega \in \so(3)$ encodes the angular velocity $\omega \in \R^3$
of the rigid-body measured in the body-fixed frame
\[
\Omega = \left( \begin{array}{ccc} 0 & -\omega_3 & \omega_2 \\
\omega_3 & 0 & -\omega_1 \\ -\omega_2 & \omega_1 & 0 \end{array} \right).
\]

In robotic applications, measurements of the attitude $R$ are
provided by exteroceptive sensor systems such as magnetometers,
accelerometers and vision systems.  The measurement of $\Omega$ is
typically obtained from on-board gyrometer systems.
\end{example}

\begin{example}\label{ex:se3}
The above example can easily be extended to model the full pose of a
rigid-body, including the position $p$ in 3D-space, by considering
invariant systems on the special Euclidean group $\SE(3)$
\cite{VikFos_2001_CDC,RehGho2003_TAC,2007_BalHamMahTru_ECC,2007_Vasconcelos_cdc}.
The special Euclidean group $\SE(3)$ has a representation as a
semidirect product \cite{bourbaki} of $\SO(3)$ and $\R^3$
$$
\SE(3)=\{ (R, p) \;|\; R \in \SO(3), p\in\R^3\},
$$
where the group product is given by
$$
(R,p)(S, q) = (RS, p+Rq) .
$$

Recall the shorthand notation
$$
T_{e}L_{\begin{pmatrix} R & p \end{pmatrix}}\begin{pmatrix} \Omega & V \end{pmatrix}=
\begin{pmatrix} R & p \end{pmatrix}\begin{pmatrix} \Omega & V \end{pmatrix}
$$
with $\begin{pmatrix} R & p \end{pmatrix}\in\SE(3)$ and
$\begin{pmatrix} \Omega & V \end{pmatrix}\in\se(3)$
and consider the left invariant system
$$
\ddt\begin{pmatrix}
R(t) & p(t)
\end{pmatrix}
=
\begin{pmatrix}
  R & p
\end{pmatrix}
\begin{pmatrix}
  \Omega & V
\end{pmatrix}
=
\begin{pmatrix}
R\Omega & RV
\end{pmatrix}
$$
with $\Omega\colon \R\rightarrow\so(3)$ and $V\colon\R\rightarrow
\R^3$ admissible. Here, $V$ is the linear velocity of the rigid-body
measured in the body-fixed frame.

In robotic applications, the measurement of $p$ is provided by
exteroceptive sensor systems such as GPS (global positioning
systems), radar or vision systems.  The measurement of linear
velocity $V$ is more challenging as there are few sensor systems
that directly measure linear velocity.  Typical measurement systems
depend on fusing differentiated exteroceptive position measurements with
integration of accelerometer measurements.
\end{example}

The two examples (Ex.~\ref{ex:so3} and \ref{ex:se3}) provide an
excellent means to demonstrate the structure of the design
principles proposed in this paper and indicate the relevance of the
results to real-world engineering problems.  It is not the purpose
of this paper, however, to enter into the details of the practical
implementation of the proposed observers for a real world system and
we will not consider issues such as choice of sensor systems or
characterisation of noise any further than the brief discussion
above and a couple of remarks later in the paper.


\section{Synchrony and error functions}\label{sec:synchrony}

In this section, we introduce the general concept of
\emph{synchrony} between pairs of systems that evolve on a given Lie
group $G$ and have a common input.  Loosely, synchrony refers to an
equivalence, but not equality, of trajectories of two dynamical
systems. More formally, we will define synchrony of two systems in
terms of an \emph{error function} $E$ that quantifies the
instantaneous difference between the system trajectories.  The pair
of systems is called \emph{$E$-synchronous} if the error $E$ is
constant along trajectories.  We show that synchronous invariant
systems on Lie-groups have a particular structure and that the error
function $E$ associated with invariant synchronous systems is also
structurally constrained.

Consider a pair of systems on $G$ driven by the same admissible
input $u$
\begin{align}
\dot X &= F_X(X,u,t) , \label{eq:sync:sys:1} \\
\dot\obs &=  F_\obs(\obs,u,t)  \label{eq:sync:sys:2}
\end{align}
with $F_X, F_\obs \colon G \times \g \times \R\rightarrow TG$,
$F_X(X,u,t)\in T_{X}G$ and $F_\obs(\obs,u,t)\in T_{\obs}G$. The
solutions of \eqref{eq:sync:sys:1} and
\eqref{eq:sync:sys:2} are denoted by $X(t; X_0, u)$ and $\obs(t;
\obs_0, u)$, respectively.

\begin{definition}
Consider the systems
\eqref{eq:sync:sys:1} and \eqref{eq:sync:sys:2}. Let
$E\colon G\times G\rightarrow M$ be a smooth \emph{error function},
$M$ a smooth manifold.
The systems \eqref{eq:sync:sys:1} and \eqref{eq:sync:sys:2} are
said to be \emph{$E$-synchronous} if for all
admissible $u\colon \R\rightarrow \g$, all initial values $X_0,
\obs_0\in G$ of \eqref{eq:sync:sys:1} and \eqref{eq:sync:sys:2}
and all $t\in \R$
$$
\ddt E(\obs(t;\obs_0,u),X(t;X_0,u)) = 0.
$$
\end{definition}

Two particularly simple error functions on a Lie
group $G$ are the canonical \emph{right invariant} error
\begin{equation}\label{eq:E_l}
   E_r(\obs,X) := \obs X^{-1}
\end{equation}
and the canonical \emph{left invariant} error
\begin{equation}\label{eq:E_r}
   E_l(\obs,X) := X^{-1} \obs .
\end{equation}
Where the arguments $\obs$ and $X$ are clear from the context we
simply write $E_r$ and $E_l$.  The label ``invariant'' refers to
simultaneous state space transformations of both systems.  That
is, for all $X,\obs,S\in G$, one has $E_r(\obs S,XS)=E_r(\obs,X)$
and an analogous result for $E_l$.

\begin{remark}
Observe that both $E_{r}$ and $E_{l}$ are non-degenerate in the sense
that the partial maps $E(\obs,\cdot)\colon G\rightarrow M$ and
$E(\cdot,X)\colon G\rightarrow M$ are global diffeomorphisms (here,
$M=G$).
\end{remark}

The next proposition shows
that any error $E$ for which two \emph{invariant} systems are
$E$-synchronous, can be factored into one of the canonical
invariant errors concatenated with a map from the Lie group into a
manifold.  Thus, the right and left invariant errors can be
thought of as the fundamental error functions for pairs of invariant
systems.

\begin{proposition} \label{prop:univ:err}
Consider the pair of invariant systems on $G$
\begin{align}
\dot X &= X u, \label{eq:syncsys_a}\\
\dot\obs &= \obs u \label{eq:syncsys_b}
\end{align}
for a single admissible input $u\colon \R\rightarrow \g$. Let $M$
be a smooth manifold. If there exists a smooth
error function $E\colon G\times G\rightarrow M$ such that the
systems \eqref{eq:syncsys_a} and \eqref{eq:syncsys_b} are
$E$-synchronous, then $E$ has the form
$$
E(\obs,X) = g( \obs X^{-1}) = g( E_r(\obs,X)),
$$
where $g\colon G\rightarrow M$ is a smooth function.

The analogous result holds for right-invariant systems
\begin{align*}
\dot X &=  v X \\
\dot \obs &= v \obs .
\end{align*}
In that case one has
$$
E(\obs,X) = g(X^{-1} \obs) = g( E_l(\obs,X))
$$
\end{proposition}
\begin{proof}
Consider the pair of systems \eqref{eq:syncsys_a} and
\eqref{eq:syncsys_b} and an error function $E$ as
above. Let $X_0,\obs_0,S \in G$ be arbitrary. Choose a smooth,
bounded curve $T\colon \R\rightarrow G$ with bounded derivatives
such that $T(0)=e$ and $T(1)=S$.  Let $(X,\obs)$ denote the solution
of the pair of systems \eqref{eq:syncsys_a} and \eqref{eq:syncsys_b}
for $u(t) = T^{-1}(t)\dot T(t)$ and
$(X,\obs)(0) = (X_0,\obs_0)$. Note that $(X,\obs)(t)=(X_0
T(t),\obs_0 T(t))$. The error function $E$ is constant on the
trajectory $(X,\obs)(t)$ and therefore $E(\obs_0,X_{0}) =
E(\obs_0S,X_{0}S)$. Since $X_0,\obs_0$ and $S$ are arbitrary, the
function $E$ is right invariant under the action of $G$, i.e.
$$
E = E\circ R_S \text{ for all } S\in G.
$$
But this implies that $E(\obs,X) = E\circ R_\xin (\obs,X) =
E(\obs\xin,e)$ for all $X,\obs\in G$. Hence $E$ has the form
$E(\obs,X) = g(\obs X^{-1})$ where $g\colon G\rightarrow M$ is the
smooth function $g(Z) = E(Z,e)$. The other case is proven by an
analogous argument.
\end{proof}

Given that the invariant errors $E_l$ and $E_r$ play such a key role
in the analysis of invariant synchronous systems it is of interest
to specialise the notion of synchrony introduced above to the
canonical errors.

\begin{definition}
A pair of systems \eqref{eq:sync:sys:1} and \eqref{eq:sync:sys:2} is
termed \emph{right synchronous} if they are $E_r$-synchronous and
\emph{left synchronous} if they are $E_l$-synchronous.
\end{definition}

In the remainder of the section, we consider the structure of a
pair of synchronous systems.  We are interested in the case where
the ``first'' system is left (resp. right) invariant \eqref{eq:sys:l}
(resp.~\eqref{eq:sys:r}).  We consider a ``partner'' system that
is either left or right synchronous and derive constraints on the
structure of the partner system.  The result will provide us with
the first half of the template for the design of non-linear
observers.  The proof of the main result
(Theorem~\ref{thm:sync:l}) uses the following lemma.

\begin{lemma} \label{lem:err}
Let $X\colon\R\rightarrow G$ and $Y\colon\R\rightarrow G$ be two smooth curves.
Then
\begin{align*}
\dot{(\xin)} Y &= - T_{\xin} R_Y T_e L_\xin T_X R_{\xin}\dot X + T_Y L_{\xin}\dot Y,\\
\dot Y \xin &= T_Y R_{\xin}\dot Y - T_{\xin} L_Y T_e L_\xin T_X R_{\xin}\dot X.\\
\end{align*}
\end{lemma}
\begin{proof}
Recall that for the inverse function $\inv(X)=X^{-1}$ we have
$T_X \inv = - T_X \left(L_{\xin} R_{\xin}\right)$.
The equations for the derivatives now follow from the usual calculus rules
for the multiplication on Lie groups.
\end{proof}

\begin{theorem} \label{thm:sync:l}
Consider the left invariant system \eqref{eq:sys:l} and let a
second system be given by the general expression
\begin{equation}\label{eq:sync}
\dot\obs = F_\obs(\obs,u,t).
\end{equation}
The systems \eqref{eq:sys:l} and \eqref{eq:sync} are right
synchronous if and only if
\begin{equation}\label{eq:l_inv_synch}
F_\obs(\obs,u,t) = \obs u .
\end{equation}
The systems are left
synchronous if and only if
\begin{equation}\label{eq:r_inv_synch}
F_\obs(\obs,u,t) = \obs \Ad_{\obs^{-1}X} u = \obs \Ad_{E_l^{-1}} u.
\end{equation}
\end{theorem}

\begin{proof}
By Lemma \ref{lem:err} we have for trajectories $X(t)$ of \eqref{eq:sys:l}
and $\obs(t)$ of \eqref{eq:sync} that
\begin{align*}
\dot{E_r} &= T_\obs R_{\xin}\dot\obs - T_{\xin} L_\obs T_e L_\xin T_X R_{\xin}\dot X \\
&= T_\obs R_{\xin} F_{\obs}(\obs,u,t) \\
&\quad - T_{\xin} L_\obs T_e L_\xin T_X R_{\xin} T_e L_X u \\
&= T_\obs R_{\xin} F_{\obs}(\obs,u,t) - T_{\xin} L_\obs T_e R_{\xin} u
\end{align*}
and
\begin{align*}
\dot{E_l} &= - T_{\xin} R_\obs T_e L_\xin T_X R_{\xin}\dot X + T_\obs L_{\xin}\dot\obs\\
&= - T_{\xin} R_\obs T_e L_\xin T_X R_{\xin} T_e L_X u \\
&\quad + T_\obs L_{\xin} F_{\obs}(\obs,u,t) \\
&= - T_{\xin} R_\obs  T_e R_{\xin}  u + T_\obs L_{\xin} F_{\obs}(\obs,u,t).
\end{align*}
Hence $E_r$ is constant if and only if
$$
F_{\obs}(\obs,u,t) = T_e L_\obs u = \obs u.
$$
and $E_l$ is constant if and only if
\begin{align*}
F_{\obs}(\obs,u,t) &= T_{\xin\obs} L_X T_{\xin} R_\obs  T_e R_{\xin} u \\
&= T_e L_\obs \Ad_{\obs^{-1}X} u = \obs \Ad_{E_l^{-1}} u.
\end{align*}
\end{proof}

In the case where the ``first'' system is right invariant
the synchronous terms for the different errors are interchanged.

\begin{theorem} \label{th:sync:r}
Consider the right invariant system \eqref{eq:sys:r} and let a second system be given by
\begin{equation}\label{eq:sync:2}
\dot\obs = F_\obs(\obs,v,t).
\end{equation}
The systems \eqref{eq:sys:r} and \eqref{eq:sync:2} are left synchronous
if and only if
$$
F_\obs(\obs,v,t) = v \obs .
$$
The systems are right synchronous if and only if
$$
F_\obs(\obs,v,t) =  (\Ad_{\obs\xin} v) \obs = (\Ad_{E_r} v) \obs.
$$
\end{theorem}


It is interesting to observe that for a left invariant system
\eqref{eq:sys:l} and choosing the right synchronous constraint,
one obtains left invariant dynamics for the partner system
\eqref{eq:l_inv_synch}.  The dynamics \eqref{eq:l_inv_synch} are
highly desirable as the fundamental component of a non-linear
observer design since they are independent of the observed system's
state $X$ and yet ensure that the observer dynamics ``track'' the
system state.  This contrasts to choosing the left synchronous
constraint, where the resulting partner system structure
\eqref{eq:r_inv_synch} depends on the error $E_l = X^{-1} \obs$
that requires knowledge of the observed system's state $X$ and cannot be
implemented in a real observer system.  Thus, it is natural to
consider right invariant errors and right synchronicity when
designing observers for left invariant systems \eqref{eq:sys:l};
and vice-versa for right invariant systems \eqref{eq:sys:r}.

\section{Internal models and innovation terms}\label{sec:innovation}
%

In this section, we define the concept of internal model
specialised for Lie group systems.  This definition is used to
show that any observer of an invariant system on a Lie group
containing an internal model, can be split into a synchronous
term, providing the internal model properties of the observer, and
a second term that we identify as an innovation term.  Thus, a
natural approach to observer design is to choose the observer as a
sum of a synchronous (internal model) term plus an innovation term,
analogous to the ``internal model plus innovation term'' design
paradigm for linear systems.

%

Consider an observed system
\begin{equation}
  \label{eq:model:sys:1}
  \dot X = F_X(X,u,t)
\end{equation}
with $F_X\colon G\times \g \times \R\rightarrow TG$, $F_{X}(X,u,t)\in T_{X}G$ and an observer
\begin{equation}
   \label{eq:model:sys:2}
   \dot\obs = F_\obs(\obs,Y,w,t)
\end{equation}
with $F_\obs\colon G\times G\times \g \times \R\rightarrow TG$,
$F_{\obs}(\obs,Y,w,t)\in T_{\obs}G$.
Note that the observer~\eqref{eq:model:sys:2} has two inputs $Y$
and $w$, with $Y$ to be fed with measurements of $X$ and $w$ to be
fed with measurements of $u$, respectively.  The idea behind an internal
model is that the observer should be able to replicate the exact
trajectory of the observed system if it is provided with exact input
information; that is, the actual initial condition of the observed
system's state and
exact measurements.  There is no a priori requirement that the system and
the observer are identical dynamical systems, only that they correspond
along certain very specific trajectories.

\begin{definition}
Consider the pair of systems \eqref{eq:model:sys:1} and
\eqref{eq:model:sys:2}. One says that \eqref{eq:model:sys:2} has
an \emph{internal model} of \eqref{eq:model:sys:1} if
for all admissible $u\colon \R\rightarrow G$, $X_0\in G$ and all
$t\in \R$
\begin{equation}
\label{eq:int:model}
\obs(t;X_0,X(t;X_0,u),u) = X(t;X_0,u) ,
\end{equation}
where $X(t;X_{0},u)$ and $\obs(t;\obs_{0},Y,w)$ denote the solutions of \eqref{eq:model:sys:1}
and \eqref{eq:model:sys:2}, respectively.
\end{definition}

Note that by Theorem~\ref{thm:sync:l} a right synchronous system
for the left invariant system~\eqref{eq:sys:l} has the form
\begin{equation}
  \label{eq:1}
  \dot\obs = \obs w
\end{equation}
with $w=u$.
It is straightforward to verify that \eqref{eq:1} has an
internal model of the original left invariant system. The
following result shows that any observer of a left invariant
system containing an internal model can be split into two parts;
the first or which is right synchronous with the original system,
providing the internal model properties of the observer, and a
second term that will be identified as an innovation term.

%

\begin{theorem}\label{th:split}
  Consider the left invariant system \eqref{eq:sys:l} and the
  observer~\eqref{eq:model:sys:2}. Assume that the observer has an
  internal model of the system. Then the right hand side of the
  observer can be written
  \begin{equation}\label{eq:obs_structure}
    F_\obs(\obs,Y,w,t)=\obs w+\alpha(\obs,Y,w,t),
  \end{equation}
  where $\alpha : G \times G \times \g \times \R \rightarrow T G$
  is a smooth function satisfying $\alpha(\obs,Y,w,t)\in T_{\obs}G$ and
  \begin{equation}\label{eq:alpha_cond}
    \alpha\left(\obs(t;X_0,X(t;X_0,u),u), X(t;X_0,u), u, t\right)=0
  \end{equation}
  for all admissible $u\colon\R\rightarrow\g$, $X_0\in G$ and $t\in\R$.
\end{theorem}

\begin{proof}
  Define $\alpha(\obs,Y,w,t)= F_\obs(\obs,Y,w,t)-\obs w$ for
  all $\obs,Y\in G$, $w\in\g$ and $t\in\R$.
  Differentiating the internal model equation~\eqref{eq:int:model} and using the
  system and observer equations yields
  \begin{multline*}
    F_\obs\left(\obs(t;X_{0},X(t;X_{0},u),u),X(t;X_{0},u),u,t\right)=\\
    X(t;X_{0},u)u
  \end{multline*}
  for all admissible $u\colon\R\rightarrow\g$, $X_0\in G$ and $t\in\R$.
  We hence get
  \begin{multline*}
    \alpha\left(\obs(t;X_0,X(t;X_0,u),u), X(t;X_0,u), u, t\right)=\\
    F_\obs\left(\obs(t;X_{0},X(t;X_{0},u),u),X(t;X_{0},u),u,t\right)-\\
    \obs(t;X_0,X(t;X_0,u),u)u=\\
    X(t;X_{0},u)u-\obs(t;X_0,X(t;X_0,u),u)u=0
  \end{multline*}
  where the last equality follows again from the definition of an internal model.
\end{proof}

A similar result can obviously be obtained for right invariant
observed systems.

The decomposition \eqref{eq:obs_structure} is analogous to the
internal model plus innovation term decomposition of linear observers.
The concept of right (resp.~left) synchronicity provides a
structural characterisation of internal models for left
(resp.~right) invariant systems.  The $\alpha$ term in Equation
\eqref{eq:obs_structure} is analogous to the innovation term in
classical linear observer design. Condition \eqref{eq:alpha_cond}
means that $\alpha$ is zero along corresponding trajectories of the
system and the observer.

\begin{definition}\label{def:innovation}
Consider the pair of systems \eqref{eq:model:sys:1} and
\eqref{eq:model:sys:2} and assume that \eqref{eq:model:sys:2} has
an internal model of \eqref{eq:model:sys:1}.
We call a map $\alpha\colon G\times G\times\g\times\R\rightarrow
TG$ an \emph{innovation term} if
\begin{itemize}
\item[\textup{(C1)}] $\alpha(\obs,Y,w,t)\in T_\obs G$ for all
  $\obs,Y\in G$, $w\in\g$, $t\in\R$ and
\item[\textup{(C2)}]
$\alpha(\obs(t;X_0,X(t;X_0,u),u), X(t;X_0,u), u, t)=0$
for all admissible $u\colon\R\rightarrow\g$, $X_{0}\in G$ and
$t\in\R$.
\end{itemize}
\end{definition}

Note that the conditions specified in Definition
\ref{def:innovation} are the least restrictive possible. In general,
an innovation term must be chosen carefully to ensure that the
trajectory of the observed system is an asymptotically stable limit set of
the observer trajectory.  We discuss design of the innovation term
in Section \ref{sec:gradients}.

Note further that Condition (C2) in
Definition~\ref{def:innovation} is in particular implied by the following
stronger condition
\begin{itemize}
\item[\textup(C2')] $\alpha(Y,Y,w,t)=0$ for all $Y\in G$,
  $w\in\g$ and $t\in\R$.
\end{itemize}
Conditions (C2) and (C2') are in general not equivalent since
there is no guarantee that for a given admissible input $u$ the trajectories of
the system (and/or those of the observer) will cover all of $G$ at any given
time.

In summary, we propose the following structure for the design of
non-linear observers for invariant systems on Lie groups.
\begin{equation}
  \label{eq:observers_l}
  \dot\obs =  \obs w_{l}+\alpha(\obs,Y,w_{l},t) \quad\text{(left observer)}
\end{equation}
This observer is intended for left invariant
observed systems of the form~\eqref{eq:sys:l}, with $Y$ receiving
measurements of $X$ and $w_{l}$ receiving measurements of $u$. Note
that this observer is in general not left
invariant, since $\alpha$ need not be left invariant.
\begin{equation}
  \label{eq:observers_r}
  \dot\obs =  w_{r}\obs+\alpha(\obs,Y,w_{r},t) \quad\text{(right observer)}
\end{equation}
Correspondingly, this observer is intended for right invariant
observed systems of the form~\eqref{eq:sys:r}, with $Y$ receiving
measurements of $X$ and $w_{r}$ receiving measurements of $v$. Note
that this observer is in general not right invariant.

%

\section{Gradient observers}\label{sec:gradients}

In this section, we present an approach to choosing innovation terms
for non-linear observers of the form \eqref{eq:observers_l}  and
\eqref{eq:observers_r} based on using the gradient
descent direction of a suitable cost function.

Let $f\colon G \times G \rightarrow \R$ be a smooth, non-negative cost
function. Furthermore, let the diagonal
$\Delta =\{ (X,X) \;|\; X\in G\}$ consist of global minima of $f$.
Recall that the Riemannian gradient of $f$ with respect to the
product metric $\langle\cdot,\cdot \rangle_p$ on $G\times G$ is
defined by
$$
\langle \grad f(\obs,Y), (\eta,\zeta) \rangle_p = df(\obs,Y)(\eta,\zeta)
$$
for all $\obs,Y\in G$, $\eta\in T_\obs G$, $\zeta\in T_Y G$. Since
we use the product metric, the gradient splits into the gradients
with respect to the first and second parameter, i.e.
\begin{multline*}
\langle \grad f(\obs,Y), (\eta,\zeta) \rangle_p = \\
\langle \grad_1 f(\obs, Y), \eta \rangle + \langle \grad_2 f(\obs,Y), \zeta
\rangle.
\end{multline*}

We propose to use the gradient of $f$ with respect to the first
parameter as an innovation term $\alpha$ in the design of
observers.  Given the observed system \eqref{eq:sys:l}
(resp.~\eqref{eq:sys:r}) then the design of the observer is
\begin{equation}\label{eq:filter:l}
\dot\obs =  \obs w_{l} - \grad_1 f(\obs,Y), \quad\text{(left observer)}
\end{equation}
respectively
\begin{equation}\label{eq:filter:r}
\dot\obs =  w_{r}\obs - \grad_1 f(\obs,Y). \quad\text{(right observer)}
\end{equation}

It remains to analyze the dynamics associated with this choice and
show that the observer trajectory asymptotically converges to the
observed system's trajectory.

\subsection{Error dynamics}

To understand the stability of the proposed observers
\eqref{eq:filter:l} and \eqref{eq:filter:r} we analyze the case
where exact measurements of the input $u$ (resp.~$v$) of
system~\eqref{eq:sys:l} (resp.~\eqref{eq:sys:r}) as well as exact
measurements of the state $X$ are available.  In terms of the
variables introduced in Section \ref{sec:innovation} one has $Y=X$
and $w_{l}=u$ (resp.~$w_{r} = v$) in the case of the observer for a left
(resp.~right) invariant system.

The following result is used in the development later in the section.

\begin{lemma} \label{lem:grad}
Let $G$ be a Lie group.  Let $f : G \times G \rightarrow \R$ be a
left invariant cost function and take a left invariant Riemmanian
metric.  Then for all $X,Y,Z\in G$
$$
T_{X} L_{Y^{-1}} \grad_1 f(X, Y Z ) =  \grad_1 f(Y^{-1}X,Z)
$$

If $f$ and the Riemmanian metric are right invariant, then
for all $X,Y,Z\in G$
$$
T_X R_{Z^{-1}} \grad_1 f(X, Y Z ) =  \grad_1 f(XZ^{-1},Y)
$$
\end{lemma}
\begin{proof}
If $f$ is left invariant, then $f\circ L_Y = f$.
Using this fact and the standard rules for transformations of Riemannian gradients, we
have that
\begin{align*}
\grad_1 f(Y^{-1}X,Z) &= \grad_1 (f\circ L_Y)(Y^{-1}X,Z )\\
&= (T_{Y^{-1}X} L_Y)^\ast \grad_1 f(X, YZ),
\end{align*}
where $(T_{Y^{-1}X} L_Y)^\ast$ denotes the Hilbert space adjoint of the linear map
$T_{Y^{-1}X} L_Y$.
Since the Riemannian metric is left invariant, we have that
$(T_{Y^{-1}X} L_Y)^\ast = T_{X} L_{Y^{-1}}$.
The right invariant case follows from an analogous argument.
\end{proof}

To analyze the asymptotic stability of the observer trajectory to
the observed system's trajectory it is convenient to consider the dynamics
of the error functions $E_r$ and $E_l$ and prove that their
trajectories converge to the identity element of the Lie group. Since
$E_{r}$ and $E_{l}$ are non-degenerate, this will imply the desired
asymptotic stability.
In particular, under suitable invariance conditions on the cost
function and the Riemannian metric the error dynamics for the
``correct'' invariant error are autonomous.

\begin{theorem} \label{thm:err:sym}
Consider the left invariant system \eqref{eq:sys:l}.  Let $f : G
\times G \rightarrow \R$ be a right invariant cost function and
take a right invariant Riemannian metric on $G$.  Consider the
left observer dynamics \eqref{eq:filter:l}.  Then the error dynamics
of the right invariant error $E_r$ \eqref{eq:E_l} are given by
\[
\dot{E_r}  = - \grad_1 f(E_r,e) .
\]

If the right invariant system \eqref{eq:sys:r} is considered with a
left invariant cost function and Riemannian metric along with the
right observer dynamics \eqref{eq:filter:r}, then the error
dynamics of the left invariant error $E_l$ \eqref{eq:E_r} are given by
$$
\dot{E_l}  = - \grad_1 f(E_l,e) .
$$
\end{theorem}

\begin{proof}
Directly from the system equations and the equations of the
left observer it follows that
\begin{align*}
\dot{E_r} &= -T_\obs R_{X^{-1}} \grad_1 f(\obs,Y)\\
&= -T_\obs R_{X^{-1}} \grad_1 f(\obs, X).
\end{align*}
Using Lemma~\ref{lem:grad} and the right invariance of $f$ and the Riemannian metric,
we get that
$$
\dot{E_r} = -\grad_1 f(\obs\xin, e) = -\grad_1 f(E_r, e) .
$$
The result for the right observer follows from an analogous argument.
\end{proof}

The gradient dynamics of the error yield the following convergence
result in the noise-free case.

\begin{theorem} \label{thm:conv}
Assume that $Y\mapsto f(Y,e)$ is a Morse-Bott function with a
global minimum at $e$ and no other local minima. If $f$ and the
Riemannian metric on $G$ are both right invariant, then both
errors $E_r$ and $E_l$ of the left observer \eqref{eq:filter:l}
converge to $e$ for generic initial conditions. Furthermore,
$f(E_r,e)$ converges monotonically to $f(e,e)$ and this
convergence is locally exponential near $e$. If $f$ and the
Riemannian metric on $G$ are both left invariant, then both errors
$E_r$ and $E_l$ of the right observer \eqref{eq:filter:r} converge
to $e$ for generic initial conditions. Furthermore, $f(E_l,e)$
converges monotonically to $f(e,e)$ and this convergence is
locally exponential near $e$.
\end{theorem}

\begin{proof}
Let us consider the left observer.
The gradient structure of the error dynamics (Theorem \ref{thm:err:sym}) and the conditions
on the cost function yield the monotonic convergence of $f(E_r,e)$ and
therefore the convergence
of $E_r$ to $e$.
As $f$ is a Morse-Bott function, the convergence of $f(E_r,e)$ is exponential near $e$.
Since convergence of $E_r$ to $e$ implies the convergence of $E_l$ to
$e$, we get the convergence of the left invariant error, too.
The result for the right observer is proved by an analogous argument.
\end{proof}

\begin{remark}
Consider the left invariant system~\eqref{eq:sys:l} and noisy
measurements $w_{l}=u+\delta$ and $Y=N_{l}X$, with additive
driving noise $\delta\in\g$ and left multiplicative
state noise $N_{l}\in G$. A straightforward
calculation yields
\begin{equation*}
  \dot{E_{r}} = \Ad_{\obs}\delta E_{r}-\grad_{1} f(E_{r},N_{l})
\end{equation*}
for the canonical right invariant error of the left
observer~\eqref{eq:filter:l}.  There is an analogous formula for
the canonical left invariant error of the right
observer~\eqref{eq:filter:r}, this time best expressed in terms of
right multiplicative state noise $Y=XN_{r}$. It is intuitively clear
that suitably bounded noise will
yield at least a practical stability result in these cases.
\end{remark}

Theorems \ref{thm:err:sym} and \ref{thm:conv} are the main results
of this section.  Together, these results provide a template for the
design of non-linear observers for invariant systems on a Lie-group,
Equations \eqref{eq:filter:l} and \eqref{eq:filter:r}.  There is
still an outstanding question of how to find suitable cost functions
that we will address in Section \ref{sub:cost_functions}. In the
remainder of this section, we consider some of the special cases
that were not addressed in Theorems \ref{thm:err:sym} and
\ref{thm:conv}.

The convergence result in Theorem \ref{thm:conv} raises the question
of what the dynamics of the other error, i.e.~$E_l$ for the left
observer and $E_r$ for the right observer, looks like and if its
cost converges monotonically to $f(e,e)$. Under suitable invariance
conditions, one obtains  the following result on the dynamics of the other
errors for each observer.

\begin{theorem} \label{thm:err:skew}
Consider the left invariant system \eqref{eq:sys:l}.  Let $f : G \times G
\rightarrow \R$ be a left invariant cost function and take a left
invariant Riemannian metric on $G$.  Consider the left observer
dynamics \eqref{eq:filter:l}.  Then the error dynamics of the
left invariant error $E_l$ are
\begin{align*}
\dot{E_l}  &= (T_e L_{E_l} - T_e R_{E_l})u - \grad_1 f(E_l,e)\\
&= (E_{l}u-uE_{l}) - \grad_1 f(E_l,e).
\end{align*}

If the right invariant system \eqref{eq:sys:r} is considered with
a right invariant cost function and Riemannian metric along with the
right observer dynamics \eqref{eq:filter:r}, then the error
dynamics of the right invariant error $E_r$ are
\begin{align*}
\dot{E_r}  &= (T_e R_{E_r} - T_e L_{E_r})v -  \grad_1 f(E_r,e)\\
&= (vE_{r}-E_{r}v) - \grad_1 f(E_r,e).
\end{align*}
\end{theorem}

\begin{proof}
Let us consider the left observer. Using Lemmas \ref{lem:err} and
\ref{lem:grad} one obtains
\begin{align*}
\dot{E_l} &=
- T_{\xin} R_\obs T_e L_\xin T_X R_{\xin} T_e L_X u
+ T_\obs L_{\xin} T_e L_\obs u \\
&\quad - T_\obs L_{\xin} \grad_1 f(\obs,Y) \\
&= (T_\obs L_{\xin} T_e L_\obs - T_{\xin} R_\obs T_e R_{\xin})u\\
&\quad - T_\obs L_{\xin} \grad_1 f(\obs,X) \\
&= (T_e L_{E_l} - T_e R_{E_l})u
- \grad_1 f(E_l,e).
\end{align*}
The statement for the right observer follows from an analogous calculation.
\end{proof}
%

It is interesting to consider the case of a bi-invariant cost
function $f$.  In this case it is possible to show that the
non-gradient term in the error dynamics
(Theorem~\ref{thm:err:skew}) is {\it passive} with respect to the
cost, i.e. an element of the kernel of the differential of the
cost.

\begin{lemma} \label{prop:passive}
Consider a Lie-group $G$ with Lie-algebra $\g$. Assume that there
exists a bi-invariant cost $f : G \times G \rightarrow \R$. Then for
all $u,v\in\g$ and any $E_l, E_r \in G$
\begin{align*}
\langle (T_e L_{E_l} - T_e R_{E_l})u ,  \grad_1 f(E_l,e)\rangle &= 0 \\
\langle (T_e R_{E_r} - T_e L_{E_r})v ,  \grad_1 f(E_r,e)\rangle &= 0
\end{align*}
\end{lemma}

\begin{proof}
For all $u\in\g$ we have by the bi-invariance of $f$
and standard transformation rules for the gradient that
\begin{align*}
\langle T_e L_{E_l}u,  \grad_1 f(E_l,e)\rangle
&=
\langle u , (T_e L_{E_l})^\ast \grad_1 f(E_l,e) \rangle \\
&=
\langle u ,  \grad_1 f(e,E_l^{-1}) \rangle \\
&=
\langle u , (T_e R_{E_l})^\ast \grad_1 f(E_l,e)\rangle \\
&= \langle T_e R_{E_l}u ,  \grad_1 f(E_l,e)\rangle .
\end{align*}
This yields the first equation.
The other equation follows from an analogous argument.
\end{proof}

\begin{proposition}\label{prop:skew_convergence}
Consider the left invariant system \eqref{eq:sys:l}.  Let $f : G \times G
\rightarrow \R$ be a bi-invariant cost function and take a
left invariant Riemannian metric on $G$.  Consider the left
observer dynamics \eqref{eq:filter:l}.  Then the time derivative
of the cost function along trajectories of \eqref{eq:sys:l} and
\eqref{eq:filter:l} is given by
\[
\ddt{f} = -|| \grad_1 f(E_l,e)||^2.
\]

The analogous result holds for a right invariant Riemannian
metric on $G$ and the right observer dynamics \eqref{eq:filter:r}.
\end{proposition}

\begin{proof}
Consider the left invariant case. The differential of $f$ is given
by
\begin{align*}
\ddt{f} & = \langle \grad_1 f(E_l,e), (T_e L_{E_l} - T_e R_{E_l})u - \grad_1 f(E_l,e) \rangle \\
& = -\langle \grad_1 f(E_l,e),  \grad_1 f(E_l,e) \rangle\\
&= -|| \grad_1 f(E_l,e)||^2
\end{align*}
by applying Lemma \ref{prop:passive}.  The analogous result for
right invariant dynamics is similarly direct.
\end{proof}

The above result can be developed into a convergence result
analogous to Theorem \ref{thm:conv}.  In practice, the case that is
of most interest is when both the cost function and the metric are
bi-invariant, a situation that occurs for the attitude estimation
example on $SO(3)$.  In this case the following convergence result
is obtained.

\begin{corollary}
Consider the left invariant system \eqref{eq:sys:l} and the right invariant system
\eqref{eq:sys:r}.  Assume that $G$ admits a bi-invariant metric
and let  $f\colon G \times G\rightarrow \R$ be a bi-invariant cost function.
Assume that $Y\mapsto f(Y,e)$ is a Morse-Bott
function with a global minimum at $e$ and no other local minima.

For both the left and right observers \eqref{eq:filter:l} and
\eqref{eq:filter:r} then $f(E_r,e)$ and $f(E_l,e)$ converge
monotonically to $f(e,e)$ for generic initial conditions. The
convergence is locally exponential near $e$.
\end{corollary}

\begin{proof}
The convergence of the cost of one type error for each filter is
shown in Corollary \ref{thm:conv}. The convergence of the other
error follows from the structure of the error dynamics, cf.
Theorem \ref{thm:err:skew}, and Proposition
\ref{prop:skew_convergence} by a straightforward Lyapunov
argument.
\end{proof}

\subsection{Example: Attitude estimation on $\SO(3)$} \label{sec:sothree}
We revisit Example~\ref{ex:so3} on the special orthogonal group
$\SO(3)$.  The left invariant attitude kinematics are
\[
\dot R = R \Omega,
\]
where $R$ denotes the attitude of a coordinate frame fixed to a
rigid body in 3D-space relative to an inertial frame and $\Omega$
encodes the angular velocity measured in the body-fixed frame. The
velocity measurements are given by $w_{l}=\Omega$ and the state
measurements are given by $Y = R$.  The right invariant and left
invariant errors have the form $E_r =\hat R R^\top$ and $E_l =
R^\top\hat R$, respectively.

We define the cost function $f(\hat R,Y)=\frac{k}{2} \| \hat R -
Y\|_F^2$, with $\|\cdot\|_F$ the Frobenius matrix norm and $k$ a
positive constant.  Since the Frobenius norm is invariant under
orthogonal transformations it follows that $f$ is bi-invariant.
Moreover, the standard Riemannian metric on $SO(3)$, induced by the
Euclidean inner product on $\R^{3\times 3}$ restricted to the
Lie-algebra of skew symmetric matrices, is also bi-invariant.  Let
us recall the well-known methods to calculate the gradient of $\hat
R\mapsto f(\hat R,Y)$, see e.g.~\cite{hm}. As we use the induced
metric on $\SO(3)$, the gradient is given by the orthogonal
projection of the Euclidean gradient in $\R^{3\times 3}$ to the
tangent space of $\SO(3)$. This projection to $T_{\hat R} SO(3)$ can
be readily calculated to be $Z \mapsto \hat R
\PP({\hat R}^{\top} Z )$, where
$\PP$ denotes the orthogonal projection onto the skew-symmetric
matrices, that is $\PP(Z)=\frac{1}{2}(Z-Z^\top)$. Since
the Euclidean gradient of $\hat R\mapsto f(\hat R,Y)$ is just
$k(\hat R-Y)$, this yields that $\grad_1 f(\hat R,Y) = k\hat R
\PP(I - {\hat R}^\top Y) = -k\hat R
\PP({\hat R}^\top Y).$ Hence we have the observer
\begin{equation*}
  \dot{\hat R} = \hat R w_{l} + k\hat R \PP({\hat R}^\top Y)
  \quad\text{(left observer)}
\end{equation*}
which coincides with the passive filter
\begin{equation*}
  \dot{\hat R} = \hat R\Omega + k\hat R \PP(\hat R^{\top}R)
\end{equation*}
proposed in \cite{2006_MahHamPfl_TAC,2005_MahHamPfl-C64} (recall that
$w_{l}=\Omega$ and $Y=R$ in the noise free case).

Starting with right invariant attitude kinematics
\begin{equation*}
  \dot R = \Gamma R
\end{equation*}
and measurements $w_{r}=\Gamma$ and $Y=R$, we have the observer
\begin{equation*}
  \dot{\hat R} = w_{r} \hat R + k\hat R \PP({\hat R}^\top Y)
  \quad\text{(right observer)}
\end{equation*}
which, using $\Gamma=\Ad_{R}\Omega$, coincides with the direct filter
\begin{equation*}
  \dot{\hat R} = \Ad_{R}\Omega\hat R + k\Ad_{\hat R}
  \PP(\hat R^{\top}R)\hat R
\end{equation*}
discussed in \cite{2006_MahHamPfl_TAC,2005_MahHamPfl-C64} but also extensively
studied over the last ten years by a range of authors
\cite{1999_Nijmeijer,ThiSan2003_TAC,Bonnabel2006_acc,2007_Tayebi_cdc,2007_Martin_cdc}.

\subsection{Construction of invariant cost functions}\label{sub:cost_functions}

In this section we investigate the question of finding invariant
cost functions on Lie groups.

For a left invariant system, Theorem \ref{thm:err:sym} shows that
the canonical right invariant error for the left observer has
gradient dynamics if the cost and Riemannian metric are right
invariant. A right invariant Riemannian metric can be easily
constructed on any Lie group by transporting a scalar product on the
Lie algebra to other tangent spaces by right translation. However,
it is a priori not clear how to obtain a right invariant cost
function, in particular if the group is non-compact.  For example,
the most natural cost function on $SE(3)$ (cf. Example~\ref{ex:se3})
would be
\begin{equation}\label{eq:simplecost:se3}
f((R,p),(Y,y))=\|R-Y\|^2+\|p-y\|^2,
\end{equation}
however, it is easily verified that this cost function is not right
invariant.  It should be noted that finding a Morse-Bott cost
function is usually fairly straightforward,
\eqref{eq:simplecost:se3} is an example of such a function, the
challenge lies in ensuring the function chosen has the desired
invariance properties.  The next proposition provides a method to
obtain an invariant cost function on a Lie group given that a
suitable Morse-Bott function on $G$ is available.

\begin{proposition}\label{prop:cost:constr}
Let $G$ be a Lie group and let $g\colon G\rightarrow \R$ be a smooth function.
Then $f\colon G\times G\rightarrow \R$, $f(X,Y)=g(XY^{-1})$ is a smooth,
right invariant function.
Furthermore, if $g$ is a Morse-Bott
function with a unique global minimum at $e$ and no further local minima,
then $Y\mapsto f(Y,e)$ is Morse-Bott
with a global minimum at $e$ and no further local minima.
\end{proposition}
\begin{proof}
The smoothness of $f$ is obvious.
For all $X,Y,Z\in G$ we have $f(XZ,YZ)=g(XY^{-1})=f(X,Y)$ and thus $f$ is right
invariant.
The second statement follows directly from $f(Y,e)=g(Y)$ for all $Y\in G$.
\end{proof}

An analogous construction yields a left-invariant cost function.
Furthermore, we can obtain a right invariant cost function from any
left invariant cost function and vice-versa from the following result.

\begin{proposition}\label{prop:inv:trick}
Let $G$ be a Lie group and let $f\colon G\times G\rightarrow\R$ be a
left invariant function.
Then $\tilde{f}\colon G\times G\rightarrow\R$ defined by
$$
\tilde{f}(X,Y)=f(X^{-1},Y^{-1})
$$
is a right invariant function.
Furthermore, if $Y\mapsto f(Y,e)$ is a Morse-Bott function such that
$e$ is the only local minimum and a global minimum, then $Y\mapsto \tilde{f}(Y,e)$
has the same properties.
\end{proposition}
\begin{proof}
The right invariance is checked by a straight forward calculation.
The second statement follows from the fact that the map
$\text{inv:}\ G\rightarrow G$, $X\mapsto X^{-1}$
is a global diffeomorphism.
\end{proof}

\subsection{Example: Pose estimation on $\SE(3)$}

As an application of the cost function construction above, let us recall
Example \ref{ex:se3}, the special Euclidean group $\SE(3)$.
The system on $\SE(3)$ is given by
$$
\ddt\begin{pmatrix}R(t) & p(t)\end{pmatrix}
=
\begin{pmatrix}
R \Omega & R V
\end{pmatrix},
$$
where $R$ resp. $p$ are the attitude resp. position of a coordinate frame
fixed to a rigid body in 3D-space relative to an inertial frame,
$\Omega$ denotes the angular velocity and $V$ denotes the linear
velocity measured in the body-fixed frame.
As mentioned before, the natural cost function \eqref{eq:simplecost:se3}
on $\SE(3)$
is not right invariant, and hence Theorem \ref{thm:conv}
cannot be applied.
However, Proposition \ref{prop:cost:constr} yields a construction procedure
for a right invariant cost function.
For this we need to choose a suitable function $g$ on $\SE(3)$.
We use
$$
g(R,p) = \frac{1}{2}\left(\|R-I\|^2 + \|p\|^2 \right).
$$
It is straightforward (although tedious) to verify that $g$ is a
Morse-Bott function with a unique global minimum at $(I,0)$ and no
further local minima. Proposition \ref{prop:cost:constr} yields the
right invariant cost function
\newcommand{\rh}{\hat{R}}
\newcommand{\ph}{\hat{p}}
\begin{align*}
f((\rh,\ph),(Y,y))
&= g((\rh,\ph)(Y,y)^{-1}) \\
&= \frac{1}{2}\left(\|\rh - Y\|^2 + \|\ph-\rh Y^\top y\|^2 \right) .
\end{align*}
Note that
\begin{equation*}
  f((\rh,\ph),(Y,y)) =
  \frac{1}{2}\left(\|\rh^\top - Y^\top\|^2 + \|-\rh^\top\ph-(-Y^\top y)\|^2 \right).
\end{equation*}
Hence we can view $f$ also as the result of applying Proposition \ref{prop:inv:trick}
to the cost function \eqref{eq:simplecost:se3}.
In order to construct a left observer, such that Theorem \ref{thm:conv}
can be applied, we also need a right invariant Riemannian metric.
Note that $T_{(I,0)} R_{(\rh,\ph)}(\Omega,V)=(\Omega \rh,V+\Omega\ph)$
for all $(\rh,\ph)\in \SE(3)$ and $(\Omega,V)\in \se(3)$.
Hence we can define a Riemannian metric by
$$
\langle (\Omega_1 \rh, V_1+\Omega_1\ph),(\Omega_2 \rh, V_2+\Omega_2\ph) \rangle
= \tr(\Omega_1^\top \Omega_2)+ V_1^\top V_2^\top
$$
for all $(\Omega_1 \rh, V_1+\Omega_1\ph), (\Omega_2 \rh, V_2+\Omega_2\ph)\in T_{(\rh,\ph)} \SE(3)$.
Note, that we have used the representation of tangent vectors by right translation of the Lie algebra
in this definition.
For our filter we have to calculate the gradient with respect to the Riemannian metric.
As a first step we derive a closed formula for the differential of $f$ with respect to the
first variable.
\begin{equation*}
d_1 f((\rh,\ph),(Y,y))(\Omega \rh,V+\Omega\ph)
= - \tr(\Omega \rh Y^T)
+ \langle \ph-\rh Y^\top y, V\rangle
\end{equation*}
where $\PP$ again denotes the orthogonal projection onto
the skew-symmetric matrices.
Hence we have the following formula for the gradient,
\begin{align*}
\lefteqn{\grad_1 f((\rh,\ph),(Y,y))} \\
&= T_{(I,0)} R_{(\rh,\ph)} (\PP(\rh Y^\top ),\ph-\rh Y^\top y) \\
&= T_{(I,0)} L_{(\rh,\ph)} \Ad_{(\rh,\ph)^{-1}}(\PP(\rh Y^\top ),\ph-\rh Y^\top y) .
\end{align*}
To obtain a representation of $\grad_1 f$ as left-translated elements of $\se(3)$ we check that
\begin{equation*}
\Ad_{(\rh,\ph)^{-1}}(\PP(\rh Y^\top ),\ph-\rh Y^\top y)=
(\PP( Y^\top\rh ),
\rh^{\top}\ph- Y^\top y + \rh^\top\PP(\rh Y^\top )\ph)
\end{equation*}
Assume now that we measure the angular and linear velocities
$(w_{\Omega},w_{V})=(\Omega,V)$ and the system state $(Y,y)=(R,p)$.
Using the construction above we get the left filter
\begin{align*}
\dot\rh &= \rh w_{\Omega} - \rh \PP( Y^\top\rh)\\
\dot\ph &= \rh w_{V} - (\ph - \rh Y^\top y) - \PP(\rh Y^\top )\ph .
\end{align*}
By Proposition \ref{prop:cost:constr}, and since $g$ was a suitable
Morse-Bott function, we can apply
Theorem \ref{thm:conv} and see that the right invariant error for this observer
converges to the identity for generic initial values.

To the best of our knowledge, the above observer on $\SE(3)$ has not
been proposed in the literature before. Previous work by two of the authors
\cite{2007_BalHamMahTru_ECC} yielded an observer on $\SE(3)$ with
exactly the same structure, but with the last term in the $\ph$-equation, i.e. 
$-\PP(\rh Y^\top)\ph$, replaced by 
$-\rh Y^{\top}\PP(\rh Y^{\top})y$. Moreover, Lyapunov-type stability
results were proven in that paper for the observers including or not including
that term, respectively.

\section{Gradient-like observers}\label{sec:grad-like}

In this section we relax the invariance requirements placed on the
cost function we have used to design the innovation term in our
gradient observers and consider general cost functions.  To provide
structure for the observer design we require that the error dynamics
exhibit autonomous gradient dynamics.  This approach leads to a
version of an \emph{internal model  principle} for invariant systems
on a Lie group.  The principle states that observers for invariant
systems that exhibit gradient dynamics for the canonical invariant
errors will contain an internal model of the observed system and can
be decomposed into a synchronous term plus an innovation term analogous
to that discussed in Section~\ref{sec:innovation}. The innovation
term is not itself a gradient term unless the cost function has the
invariance properties discussed in Section \ref{sec:gradients} and
we term the resulting observers \emph{gradient-like}.  This result
sharpens a related structural result in \cite{Bonnabel2005_ch} (see
also \cite{2008_Bonnabel_TAC}), where less stringent conditions were
placed on the error dynamics.  Sharpening these requirements
allows us to derive an almost global convergence result in a very
general context.

\begin{theorem}\label{th:gl:filter}
Let $f\colon G\times G\rightarrow\R$ be a smooth cost function.
Consider a general left observer for system~\eqref{eq:sys:l},
\begin{equation}\label{eq:abs:inno}
\dot\obs = F_\obs(\obs,Y,w_{l},t)
\end{equation}
with measurements $Y=X$ and $w_{l}=u$.  Then the canonical right invariant error $E_{r}$
displays gradient dynamics
\begin{equation}\label{eq:err:grad}
\dot{E_r} = - \grad_1 f(E_r,e),
\end{equation}
if and only if
\[
F_\obs(\obs,Y,w_{l},t) = \obs w_{l} - T_{\obs Y^{-1}}R_{Y} \grad_1 f(\obs Y^{-1},e) .
\]

The analogous result holds for right observers for
system~\eqref{eq:sys:r} with measurements $Y=X$
and $w_{r}=v$, and the canonical left invariant error $E_{l}$. In this case, the
dynamics of the left invariant error are gradient dynamics
\begin{equation}\label{eq:err:grad:l}
\dot{E_l} = - \grad_1 f(E_l,e),
\end{equation}
if and only if
\[
F_\obs(\obs,Y,w_{r},t) = w_{r}\obs - T_{Y^{-1}\obs} L_Y \grad_1 f(Y^{-1} \obs,e) .
\]
\end{theorem}

\begin{proof}
Assume that we have an observer \eqref{eq:abs:inno}
with error dynamics \eqref{eq:err:grad}.
We can split $F_\obs(\obs,Y,w_{l},t)=\obs w_{l}+\alpha$
into the synchronization term $\obs w_{l}$ and a
remainder term $\alpha$ by defining $\alpha=F_\obs(\obs,Y,w_{l},t)-\obs w_{l}$.
By Lemma \ref{lem:err} and using the same argument as in the proof of
Theorem \ref{thm:sync:l}, as well as $w_{l}=u$, we see that
$$
\dot{E_r} = T_\obs R_{X^{-1}} \alpha.
$$
We see that $\alpha = - T_{\obs Y^{-1}}R_{Y} \grad_1 f(\obs Y^{-1},e)$
by using $X=Y$ and \eqref{eq:err:grad}. In particular,
$\alpha$ is an innovation term in the sense of Section
\ref{sec:innovation}. If on the other hand $F_\obs(\obs,Y,w_{l},t)$
has the form as given in the Theorem then a straightforward
calculation shows that we get the error dynamics
\eqref{eq:err:grad}. The statement for the error dynamics
\eqref{eq:err:grad:l} follows analogously.
\end{proof}

Theorem \ref{th:gl:filter} yields the two gradient-like observers
\begin{equation}\label{eq:filter:gl:l}
  \dot\obs = \obs w_{l} - T_{\obs Y^{-1}} R_Y \grad_1 f(\obs Y^{-1},e)
\end{equation}
and
\begin{equation}\label{eq:filter:gl:r}
  \dot\obs = w_{r} \obs  - T_{Y^{-1}\obs} L_Y \grad_1 f(Y^{-1} \obs,e).
\end{equation}
Note that by Lemma \ref{lem:grad} these observers coincide with the
gradient observers \eqref{eq:filter:l} and \eqref{eq:filter:r},
respectively, if the cost function and the Riemannian metric are
right, respectively left, invariant.

\begin{corollary}\label{cor:gl:filter:conv}
Assume that $Y\mapsto f(Y,e)$ is a Morse-Bott function with a global minimum at $e$
and no other local minima.
Both errors \textup{(}$E_r$ and $E_l$\textup{)} of the left observer
\eqref{eq:filter:gl:l} and both errors of the right observer
\eqref{eq:filter:gl:r} converge to $e$ for generic initial conditions.
Furthermore, $f(E_r,e)$ converges monotonically to $f(e,e)$ for the left filter
and $f(E_l,e)$ converges monotonically to $f(e,e)$ for the right
filter, and this convergence is locally exponential near $e$ in both cases.
\end{corollary}

\begin{proof}
The proof is entirely analogous to that of Theorem \ref{thm:conv}.
\end{proof}

\begin{remark}
In the case of noisy measurements $w_{l}=u+\delta$ and $Y=XN_{r}$, a
straightforward calculation shows
$$
\dot{E_r} = \Ad_{\obs}\delta E_{r}
- T_{E_r N_r^{-1}} R_{\xin N_r X} \grad_1 f(E_r N_r^{-1} ,e)
$$
for the canonical right invariant error of  the left gradient-like
observer~\eqref{eq:filter:gl:l}.  The formula for the canonical left
invariant error of the right gradient-like
observer~\eqref{eq:filter:gl:r} is analogous.  It is intuitively clear that these
systems will have strong practical stability properties.
\end{remark}

\section{Conclusion}\label{sec:conclusions}

This paper provides a coherent theory of nonlinear observer design
for systems with left (resp.~right) invariant kinematics on a Lie
group for which one has a right (resp.~left) invariant,
non-degenerate, Morse-Bott cost function and full measurements.
The key contributions are the observer equations
\eqref{eq:filter:l} and \eqref{eq:filter:r} along with Theorems
\ref{thm:err:sym} and \ref{thm:conv}.  The results in Section
\ref{sub:cost_functions} are of practical importance in generating
invariant cost functions. Finally, the results presented in
Section \ref{sec:grad-like} provide a practical design methodology
in the case where a non-invariant cost function is considered.  A
limitation of the approach described in this paper is the
requirement for full measurement of both state and velocity.  In
work in progress, we are investigating the structure of observers
for kinematic systems on Lie groups with partial state
measurements.

\bibliographystyle{IEEEtran}

\bibliography{glfilter_4}

\end{document}